\documentclass{amsart}
\usepackage{amsmath,amsthm,amssymb,amscd}
\newtheorem{thm}{Theorem}

\newtheorem{prop}{Proposition}
\newtheorem{lem}{Lemma}

\usepackage[dvipdfmx]{graphicx}

\begin{document}
\title
[ The Banach-Tarski paradox for flag manifolds]
{ The Banach-Tarski paradox for flag manifolds}
\author{Yohei Komori}
\author{Yuriko Umemoto}
\address[Y. Komori]{Osaka City University Advanced Mathematical Institute and
Department of Mathematics, Osaka City University, 558-8585, Osaka, Japan}
\email{komori@sci.osaka-cu.ac.jp}
\address[Y. Umemoto]{Department of Mathematics, Osaka City University, 558-8585, Osaka, Japan}
\email{yuriko.ummt.77@gmail.com}
\subjclass[2010]{Primary~22F30, Secondary~14M15}
\keywords{equidecomposability, paradoxical action, flag manifold.}
\date{}
\thanks{This work is partially supported by the JSPS Institutional Program for 
Young Researcher Overseas Visits
`` Promoting international young researchers in mathematics and 
mathematical sciences led by OCAMI ".
}
\begin{abstract}
The famous Banach-Tarski paradox claims that the three dimensional rotation group $SO(3)$ acts on the two dimensional sphere $S^2$ paradoxically.
In this paper, we generalize their result to show that the  classical group $G(n,K)$ acts on the flag manifold $F(d_1, d_2, \cdots, d_k, K)$ paradoxically.
\end{abstract}
\maketitle

\section{Introduction}

Let $X$ be a non-empty set on which a group $G$ acts.
In this case $X$ is called a {\bf G-space}.
In particular if any non-identity element of $G$  acts on $X$ without fixed points, 
we say that $G$ acts on $X$ {\bf freely}.
For example  $G$ acts naturally on itself by left translation freely.

Non-empty subsets $A$ and $B$ of $G$-space $X$ are called {\bf $G$-equidecomposable} if
there exist finite elements of $G$, $g_1, g_2, \cdots, g_n \in G$,
and n-partitions of $A$ and $B$ respectively,  $A=\sqcup_{i=1}^{n}  A_i$, 
$B=\sqcup_{i=1}^{n}  B_i$  such that
$$
A_i=g_i B_i, \;\; \forall i \in \{1,2, \cdots, n \}.
$$

Non-empty subset $E$ of $G$-space $X$ is called {\bf $G$-paradoxical} if
there exist disjoint subsets $A$ and $B$ of $E$ such that
$A$ is $G$-equidecomposable to $E$ while $B$ is also $G$-equidecomposable to $E$.

For example, the rank two free group $F_2$ is $F_2$-paradoxical (as $F_2$-space):
In practice take  free generators $a$ and $b$.
For $x \in \{ a, b, a^{-1}, b^{-1} \}$, let $W(x)$ be the set of reduced words whose prefix is the letter  $x$.
Then
$((W(a)  \sqcup  W(a^{-1})) \sqcup ( (W(b)  \sqcup  W(b^{-1})) \subset  F_2$
and $F_2=W(a) \sqcup a W(a^{-1}) =  W(b) \sqcup b W(b^{-1})$.
The following claims are easy consequences from the definitions \cite{sw}.
\begin{prop}
\label{prop: key}
\begin{enumerate}
\item
Suppose that $H$ is a subgroup of $G$ and $X$ is a $G$-space.
If $X$ is $H$-paradoxical, then $X$ is also $G$-paradoxical.
\item
If $G$-spaces $X$ and $Y$ are disjoint and $G$-paradoxical, then
the disjoint union $X \sqcup Y$ is also $G$-paradoxical.
\item
For non-empty subsets $A$ and $B$ of $G$-space $X$,
suppose that $A$ and $B$ are $G$-equidecomposable and $A$ is $G$-paradoxical.
Then $B$ is also $G$-paradoxical.
\item
If $G$ is $G$-paradoxical (as $G$-space) and acts on $X$ freely,
then $X$ is also $G$-paradoxical.
\item
Suppose that $X$ and $Y$ are $G$-spaces and there is a $G$-equivariant map from $X$ to $Y$.
If $Y$ is $G$-paradoxical, then $X$ is also $G$-paradoxical.
\end{enumerate}
\end{prop}

Banach and Tarski showed that the 3-dimensional rotation group $SO(3)$ acts on the 2-dimensional sphere $S^2$ 
paradoxically, which is known as the Banach-Tarski paradox (\cite{sw} Corollary 3.10).
Since $S^2$ is the homogeneous space $SO(3)/SO(2)=O(3)/O(2)$, 
Several generalizations of this result to other homogeneous spaces were considered:
in practice it is known that the $(n-1)$-dimensional sphere $S^{n-1}=O(n)/O(n-1)$ and the real projective space $\mathbf{R}P^{n-1}=O(n)/O(1) \times O(n-1)$ are $O(n)$-paradoxical
 \cite{sw}.
In this paper we will show that for any $n \geq 3$ and any sequence of natural numbers $(n_1, n_2, \cdots, n_k)$ satisfying $n_1+n_2+\cdots +n_k=n$,
the real flag manifold $O(n)/O(n_1) \times O(n_2) \times \cdots \times O(n_k)$
is $O(n)$-paradoxical.
More generally we will consider complex and quotanionic  flag manifolds also.

We now describe the contents of this paper.
In section 2, we will review the definitions of projective spaces, Grassmann manifolds and flag manifolds 
over the real number field $\mathbb{R}$ , the complex number field $\mathbb{C}$, and the quotanion algebra  $\mathbb{H}$ as homogeneous spaces
of the classical groups $O(n)$,  $U(n)$, and $Sp(n)$ respectively.
We will show our main theorem for partial flag manifolds and explain that it reduces to the same result for Grassmann manifolds
which will be proved in the final section.
In section 3 we will review the idea of the proof of Banach-Tarski paradox for spheres following \cite{sw} which  we will use in section 4.
In section 4 we will prove  the Banach-Tarski paradox for projective spaces, and by using this
we will show  the Banach-Tarski paradox for Grassmann manifolds in section 5.

\section{Notations and the main result}
For $K=\mathbb{R, C}$, and $\mathbb{H}$, 
the $n$-dimensional right $K$-vector space $K^n$
has the following inner product:
$$
(x,y):=\bar{x}_1y_1+ \cdots + \bar{x}_ny_n
$$
for $x,y \in K$, which defines the metric $d(x,y)$ on $K^n$ by
$$
d(x,y):=\sqrt{(x-y, x-y)}.
$$
The isometry group of this metric is the compact Lie group $O(n), U(n)$, and $Sp(n)$ for $K=\mathbb{R, C}$, and $\mathbb{H}$ respectively, 
which we will denote by $G(n,K)$ for simplicity.
In this paper we will consider the paradoxical action of $G(n,K)$ on the following homogeneous spaces of $G(n,K)$
(in practice they are symmetric spaces in the sense of differential geometry).

First let us denote the set of all lines through the origin (i.e. 1-dimensional right $K$-subspaces) in $K^n$ by $KP^{n-1}$ and call it the {\bf $(n-1)$-dimensional K-projective space}.
$G(n,K)$ acts on $KP^{n-1}$ transitively so that $KP^{n-1}$ becomes the $G(n,K)$-homogeneous space as follows:
$$
KP^{n-1}=G(n,K)/G(1,K) \times G(n-1,K).
$$

Next we consider the set of all $d$-dimensional  right $K$-subspaces in  $K^n$ by $Gr_d K^n$ and call it a {\bf Grassmann manifold}.
The $(n-1)$-dimensional K-projective space $KP^{n-1}$ is the Grassmann manifold $Gr_1 K^n$.
$G(n,K)$ acts on $Gr_d K^n$ transitively so that $Gr_d K^n$ becomes the $G(n,K)$-homogeneous space as follows:
$$
Gr_d K^n=G(n,K)/G(d,K) \times G(n-d,K).
$$

Finally we define a {\bf partial flag of index $(d_1, d_2, \cdots, d_k)$} by a strictly increasing sequence of  right $K$-subspaces of $K^n$
\cite{aa, bc, hh}:
$$
\{ 0 \}=V_0 \subset V_1 \subset \cdots \subset V_k=K^n, \;\; d_i:=dim_K V_i.
$$
Let us denote the set of all partial flags of index $(d_1, d_2, \cdots, d_k)$ by $F(d_1, d_2, \cdots, d_k, \\K)$ and call it a {\bf flag manifold}.
The Grassmann manifold $Gr_d K^n$ is the flag manifold $F(d, n, K)$.
$G(n,K)$ acts on $F(d_1, d_2, \cdots, d_k, K)$ transitively so that $F(d_1, d_2, \cdots, d_k, K)$ becomes the $G(n,K)$-homogeneous space as follows:
$$
F(d_1, d_2, \cdots, d_k, K)=G(n,K)/G(n_1,K) \times G(n_2,K) \times \cdots \times G(n_k,K)
$$
where $n_i$ denotes $d_i-d_{i-1}$.
Next result is a consequence from the definitions:
\begin{prop}
\label{prop: flag}
By taking the i-th component $V_i$ of a partial flag
$V_0 \subset V_1 \subset \cdots \subset V_k$,
there is a $G(n,K)$-equivariant map from $F(d_1, d_2, \cdots, d_k, K)$ to $Gr_{d_i} K^n$.
\end{prop}

Our main purpose is to show that the Banach-Tarski paradox holds for the action of
$G(n,K)$ on $F(d_1, d_2, \cdots, d_k, K)$:

\begin{thm}
\label{thm: flag}
Let $n_K \in \mathbb{N}$ be equal to $3$ when $K=\mathbb{R}$, and equal to $2$ when $K=\mathbb{C}$ or $\mathbb{H}$.
Then for any $n \geq n_K$ and any sequence $(d_1, d_2, \cdots, d_k)$ satisfying $0<d_1<d_2<\cdots <d_k=n$,
the flag manifold $F(d_1, d_2, \cdots, d_k, K)$ is $G(n,K)$-paradoxical.
\end{thm}

By means of Proposition \ref{prop: key}.(5) and Proposition \ref{prop: flag}, Theorem \ref{thm: flag} will be a consequence of the following theorem
which we will prove in section 5:

\begin{thm}
For any $n \geq n_K$ and any $1 \leq k \leq n-1$,
$Gr_k K^n$ is $G(n, K)$-paradoxical
where $n_K=3$ when $K=\mathbb{R}$ and $n_K=2$ when $K=\mathbb{C}$ or $\mathbb{H}$.
\end{thm}

\section{Spheres}
The linear action of $SO(n+1)$ on $\mathbb{R}^{n+1}$ induces the action $SO(n+1)$ on the
n-dimensional sphere
$$
S^n:=\{ (x_1, x_2, \cdots, x_{n+1}) \in \mathbb{R}^{n+1} \; | \; \sum_{k=1}^{n+1} x_k^2=1 \}.
$$

In this section, following \cite{sw} Theorem 5.1, we will prove

\begin{thm}
$S^n$ is $SO(n+1)$-paradoxical
for all $n \geq 2$.
\end{thm}

We will show this claim by induction on $n$:
\begin{enumerate}
\item
First we consider the case $n=2$:
\begin{prop}(\cite{sw} Theorem 2.1)
$SO(3)$ contains a subgroup isomorphic to the rank 2 free group $F_2$.
More precisely the subgroup $H$ of $SO(3)$ generated by the following matrices
$A$ and $B$ is isomorphic to $F_2$:
$$
A=
\left(
\begin{array}{ccc}
\frac{1}{3}  & -\frac{2\sqrt{2}}{3}  &  0 \\
\frac{2\sqrt{2}}{3}  & \frac{1}{3}  &  0 \\
0  & 0  &   1
\end{array}
\right), \;\;\;\;
B=
\left(
\begin{array}{ccc}
1 & 0 & 0\\
0 & \frac{1}{3}  & -\frac{2\sqrt{2}}{3}\\
0 & \frac{2\sqrt{2}}{3}  & \frac{1}{3}
\end{array}
\right).
$$
\end{prop}

For any $g \in H - \{ id\}$, let $\ell_g$ be the rotation axis for the linear action of $g$ on $\mathbb{R}^3$.
Then the countable set 
$D := \bigcup_{g \in H - \{id\}} (\ell_g \cap S^2)$ is $H$-invariant so that $H$ acts on $S^2 - D$ freely.
Hence Proposition \ref{prop: key}.(4) implies
\begin{prop}(\cite{sw} Theorem 2.3)
\label{prop: para}
$S^2 - D$ is $SO(3)$-paradoxical.
\end{prop}

\begin{prop}(\cite{sw} Theorem 3.9)
\label{prop: equi}
$S^2$ and $S^2 - D$ are $SO(3)$-equidecomposable.
\end{prop}
(Proof.)\\
For any $m \neq n \in \mathbb{N}$, there exists $g \in SO(3) -H$ such that
 $g^m(D) \cap g^n(D) = \emptyset$.
 Put $A:= \cup_{n=0}^\infty g^n(D)$.
 Then $A$ is a countable set since $D$ is countable.
Hence
$S^2=(S^2-A) \sqcup A$ and $S^2-D=(S^2-A) \sqcup g(A)$ imply that
$S^2$ and $S^2 - D$ are $SO(3)$-equidecomposable.

Therefore by means of Proposition \ref{prop: key}.3,
$S^2$ is $SO(3)$-paradoxical, which is known as the Banach-Tarski paradox (\cite{sw} Corollary 3.10).

\item
By the induction hypothesis, 
we assume that there exists $k  \geq 2$ such that
$S^k$ is $SO(k+1)$-paradoxical.

Let  $H_{k+1} \subset \mathbb{R}^{k+2}$ be the image of $\mathbb{R}^{k+1}$
under the natural embedding
\begin{eqnarray*}
\mathbb{R}^{k+1} & \rightarrow &  \mathbb{R}^{k+2}\\
 (x_1, x_2, \cdots, x_{k+1}) & \mapsto & (x_1, \cdots, x_{k+1}, 0).
\end{eqnarray*}
Then
$S^k \subset \mathbb{R}^{k+1}$ can be identified with
$S^{k+1} \cap H_{k+1}$.
This	identification realizes the action of
$SO(k+1)$ on $S^k$ as the action of the subgroup
$SO(k+1)^*:=
\left(
\begin{array}{cc}
SO(k+1)  &   0_{ k+1,1}   \\
0_{1,k+1}  &   1
\end{array}
\right)$
of $SO(k+2)$
on $S^{k+1} \cap H_{k+1}$
where $0_{ k+1,1}$ and $0_{1,k+1}$ are zero matrices of sizes $(k+1) \times 1$ and
$1 \times (k+1)$ respectively.
By the induction hypothesis,
$S^{k+1} \cap H_{k+1}$ is $SO(k+1)^*$-paradoxical.

The natural projection
\begin{eqnarray*}
S^{k+1} - \{ (0, \cdots, \pm1) \} & \rightarrow &  S^{k+1} \cap H_{k+1}\\
 (x_1, \cdots, x_{k+1}, x_{k+2}) & \mapsto &  \frac{1}{\sqrt{\sum_{i=1}^{k+1} x_i^2}}(x_1, \cdots, x_{k+1}, 0)
\end{eqnarray*}
is $SO(k+1)^*$-equivariant.

\begin{figure}[h]
 \begin{center}
 \includegraphics [width=200pt, clip]{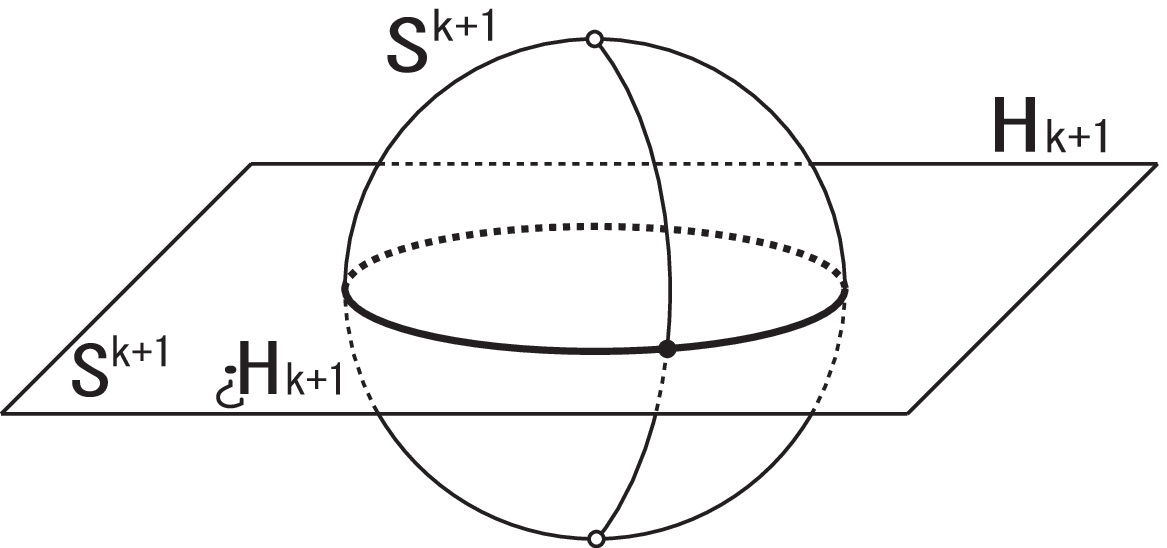}
\end{center}
\caption{}
\end{figure}

Hence Proposition \ref{prop: key}.5 implies that $S^{k+1} - \{ (0, \cdots, \pm1) \} $ is $SO(k+1)^*$-paradoxical, in particular $SO(k+2)$-paradoxical
by Proposition \ref{prop: key}.(1).
Moreover by the same argument of Proposition \ref{prop: equi}, $S^{k+1} - \{ (0, \cdots, \pm1) \} $ and 
$S^{k+1}$
are $SO(k+2)$-equidecomposable.
Therefore Proposition \ref{prop: key}.3 implies that $S^{k+1}$ is $SO(k+2)$-paradoxical.
\end{enumerate}

\section{Projective spaces}

\begin{thm}
\label{thm: proj}
Let $n_K \in \mathbb{N}$ be equal to $3$ when $K=\mathbb{R}$, and equal to $2$ when $K=\mathbb{C}$ or $\mathbb{H}$.
Then the $(n-1)$-dimensional $K$-projective space
$KP^{n-1}$ is $G(n,K)$-paradoxical for any $n \geq n_K$.
\end{thm}

In the following we will show this claim by induction on $n$:
\begin{enumerate}
\item
First we consider the case $n=n_K$.
When $K=\mathbb{R}$, then it claims that
${\mathbb R}P^2$ is $SO(3)$-paradoxical.

Any line $\ell$ in ${\mathbb R}^n$ passing through the origin
intersects the $n$-dimensional sphere at antipodal points  $\pm Q$.
Hence there exists  a natural 2 to 1 surjective map
$$
\pi: S^n \rightarrow  {\mathbb R}P^{n}
$$ 
which identifies antipodal points  $\pm Q$.
Also because
$SO(n+1)$ action on $S^n$ comes from the linear action of $SO(n+1)$ on $\mathbb{R}^{n+1}$,
for $M \in SO(n+1)$ and $Q \in S^n$,
$M(-Q)=-M(Q)$, which means that the $SO(n+1)$ action on $S^n$ induces the action on
${\mathbb R}P^{n}$  so that $\pi$ is $SO(n+1)$-equivariant.
Then by analogy with Proposition \ref{prop: para}
\begin{prop}
${\mathbb R}P^{2}-\pi(D)$
is $SO(3)$-paradoxical.
\end{prop}
Also by similar arguments of Proposition \ref{prop: equi}
\begin{prop}
${\mathbb R}P^{2}$ and ${\mathbb R}P^{2}-\pi(D)$ are $SO(3)$-equidecomposable.
\end{prop}

Therefore by means of Proposition \ref{prop: key}.(3),
${\mathbb R}P^{2}$  is $SO(3)$-paradoxical.

Next we consider the case when $K=\mathbb{C}$ and $\mathbb{H}$:
we will show that
$KP^1$ is $G(2, K)$-paradoxical.
The linear action of $G(2, K)$ on $K^2$ reduces to the action of $G(2, K)$ on $KP^1$.
By means of the stereographic projection it reduces to the action of $SO(3)$ on $S^2$
when $K=\mathbb{C}$, and the action of $SO(5)$ on $S^4$ when $K=\mathbb{H}$ \cite{pg}.
Hence from the result of section 2, we can conclude our claim.

\item
By the induction hypothesis, 
we assume that there exists $k  \geq n_K$ such that
$KP^{k-1}$ is $G(k, K)$-paradoxical.

Let  $H_{k} \subset K^{k+1}$ be the image of $K^{k}$
under the natural embedding
\begin{eqnarray*}
K^k & \rightarrow & K^{k+1}\\
 (x_1, \cdots, x_k) & \mapsto & (x_1, \cdots, x_k, 0).
\end{eqnarray*}
Then
$KP^{k-1}$ 
can be identified with
$(KP^{k-1})^*:=\{ \ell \in KP^k \; | \; \ell \subset H_{k} \}$.
This	identification realizes the action of
$G(k,K)$ on $KP^{k-1}$ as the action of
the subgroup
$G(k,K)^*=
\left(
\begin{array}{cc}
G(k,K)  &   0_{k,1}   \\
 0_{1,k} &   1
\end{array}
\right)$
of $G(k+1,K)$
on $(KP^{k-1})^*$
where $0_{ k,1}$ and $0_{1,k}$ are zero matrices of sizes $k \times 1$ and
$1 \times k$ respectively.
By the induction hypothesis,
$(KP^{k-1})^*$ is $G(k,K)^*$-paradoxical.

The natural projection
\begin{eqnarray*}
P_{k+1}: K^{k+1} & \rightarrow &   K^k\\
(x_1, \cdots, x_k, x_{k+1}) & \mapsto & (x_1, \cdots, x_k)
\end{eqnarray*}
induces the following $G(k,K)^*$-equivariant map
\begin{eqnarray*}
KP^k - \{ x_{k+1}\text{-axis} \} & \rightarrow & (KP^{k-1})^*\\
\ell & \mapsto & P_{k+1}(\ell).
\end{eqnarray*}

\begin{figure}[h]
 \begin{center}
 \includegraphics [width=170pt, clip]{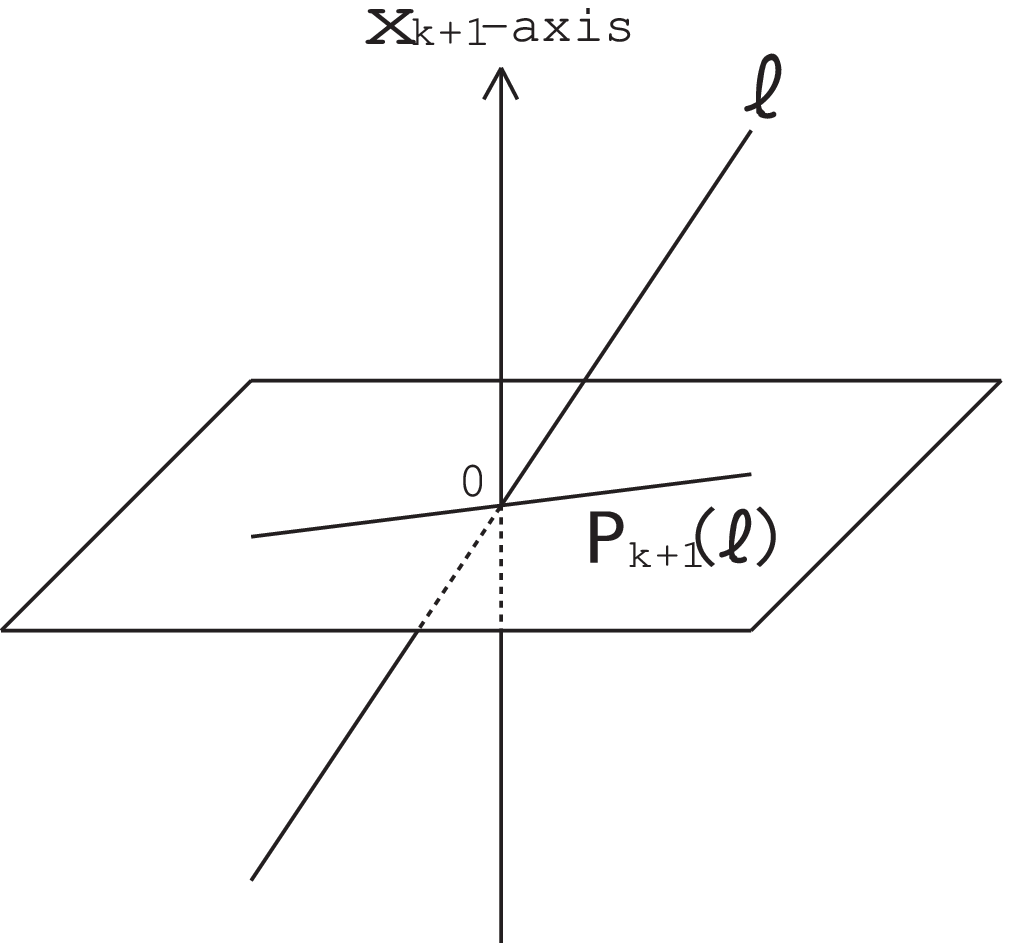}
\end{center}
\caption{}
\end{figure}

Hence Proposition \ref{prop: key}.(5) implies that $KP^k - \{ x_{k+1}\text{-axis} \}$ is $G(k,K)^*$-paradoxical, in particular $G(k+1, K)$-paradoxical
by Proposition  \ref{prop: key}.(1).
Moreover by the same argument of Proposition \ref{prop: equi}, $KP^k - \{ x_{k+1}\text{-axis} \}$ and 
$KP^k$
are $G(k+1, K)$-equidecomposable.
Therefore Proposition \ref{prop: key}.(3) implies that $KP^k$ is $G(k+1, K)$-paradoxical.
\end{enumerate}

\section{Grassmann manifolds}

In this section we will prove the Banach-Tarski paradox for Grassmann manifolds
which induces our main result Theorem 1 appeared in section 1.
Key idea is the following duality between Grassmann manifolds:
\begin{prop}
\label{prop: dual}
Let $\varphi: Gr_k K^n \rightarrow Gr_{n-k} K^n$
be the map defined by $\varphi(H)=H^{\perp}$
where  $H^{\perp}$ is the orthogonal complement of  $H$ in $K^n$.
Then
$\varphi$ is a $G(n,K)$-equivariant homeomorphism.
\end{prop}

\begin{thm}
For any $n \geq n_K$ and any $1 \leq k \leq n-1$,
$Gr_k K^n$ is $G(n, K)$-paradoxical
where $n_K=3$ when $K=\mathbb{R}$ and $n_K=2$ when $K=\mathbb{C}$ or $\mathbb{H}$.
\end{thm}

In the following we will show this claim by induction on $n$:
\begin{enumerate}
\item
First we consider the case $n=n_K$.
When $K=\mathbb{R}$, then it claims that
$Gr_k {\mathbb R}^3$ is $SO(3)$-paradoxical for $k=1,2$.
When $k=1$, $Gr_1{\mathbb R}^3=\mathbb{R}P^2$ which is $SO(3)$-paradoxical
by Theorem \ref{thm: proj}.
When $k=2$, 
because of Proposition \ref{prop: key}.(3) and Proposition \ref{prop: dual}, 
$Gr_2{\mathbb R}^3$ is also $SO(3)$-paradoxical.
When $K=\mathbb{C}$ or $\mathbb{H}$,  then it claims that
$Gr_1 K^2=KP^2$ is $G(2, K)$-paradoxical which is also proved 
by Theorem \ref{thm: proj}.
\item
By the induction hypothesis, 
we assume that there exists $n_0 >n_K$ such that
for any $n$ satisfying $n_K \leq n < n_0$ and
any $k$ satisfying $1 \leq k \leq n-1$,
$Gr_k K^n$ is $G(n, K)$-paradoxical.
Under this assumption, we will show that
for any  $k$ satisfying $1 \leq k \leq n_0-1$,
$Gr_k K^{n_0}$ is $G(n_0, K)$-paradoxical.
\begin{enumerate}
\item
For $k=1$, $Gr_1K^{n_0}=KP^{n_0-1}$ 
which is $G(n_0, K)$-paradoxical by Theorem \ref{thm: proj}.
\item
Next we consider the case $2 \leq k \leq n_0/2$.

Let  $H_{n_0+1-k} \subset K^{n_0}$ be the image of $K^{n_0+1-k}$
under the natural embedding
\begin{eqnarray*}
K^{n_0+1-k} & \rightarrow & K^{n_0}\\
 (x_1, \cdots, x_{n_0+1-k}) & \mapsto & (x_1, \cdots, x_{n_0+1-k}, 0, \cdots, 0)
\end{eqnarray*}
Then 
$Gr_k K^{n_0+1-k}$ can be identified with
$(Gr_k K^{n_0+1-k})^*:=\{ V \in Gr_k K^{n_0} \; | \; V \subset H_{n_0+1-k} \}$.
This identification realizes 
the action of $G(n_0+1-k, K)$ on $Gr_k K^{n_0+1-k}$
as the action of the subgroup
$G(n_0+1-k, K)^*:=
\left(
\begin{array}{cc}
G(n_0+1-k, K)  &   0_{ n_0+1-k, k-1}  \\
0_{k-1,n_0+1-k} &   E_{k-1,k-1}
\end{array}
\right)$
of $G(n_0, K)$
on $(Gr_k K^{n_0+1-k})^*$,
where $0_{ n_0+1-k, k-1}$ and $0_{k-1,n_0+1-k}$ are zero matrices of sizes $(n_0+1-k) \times (k-1)$ and
$(k-1)\times (n_0+1-k)$ respectively, and 
$E_{k-1,k-1}$ is the identity matrix of size $(k-1)\times (k-1)$.

Similarly by identifying 
$Gr_1 K^{n_0+1-k}$ with
$(Gr_1 K^{n_0+1-k})^*:=\{ L \in Gr_1K^{n_0} \; | \; L \subset H_{n_0+1-k} \}$,
the action of
$G(n_0+1-k, K)$ on $Gr_1 K^{n_0+1-k}$ can be realized by the action of
 $G(n_0+1-k, K)^*$ 
on $(Gr_1K^{n_0+1-k})^*$.

\begin{lem}
The following map is $G(n_0+1-k, K)^*$-equivariant.
\begin{eqnarray*}
Gr_k K^{n_0} - (Gr_k K^{n_0+1-k})^* & \rightarrow & (Gr_1 K^{n_0+1-k})^*\\
V &  \mapsto & V \cap H_{n_0+1-k}.
\end{eqnarray*}

\begin{figure}[h]
 \begin{center}
 \includegraphics [width=170pt, clip]{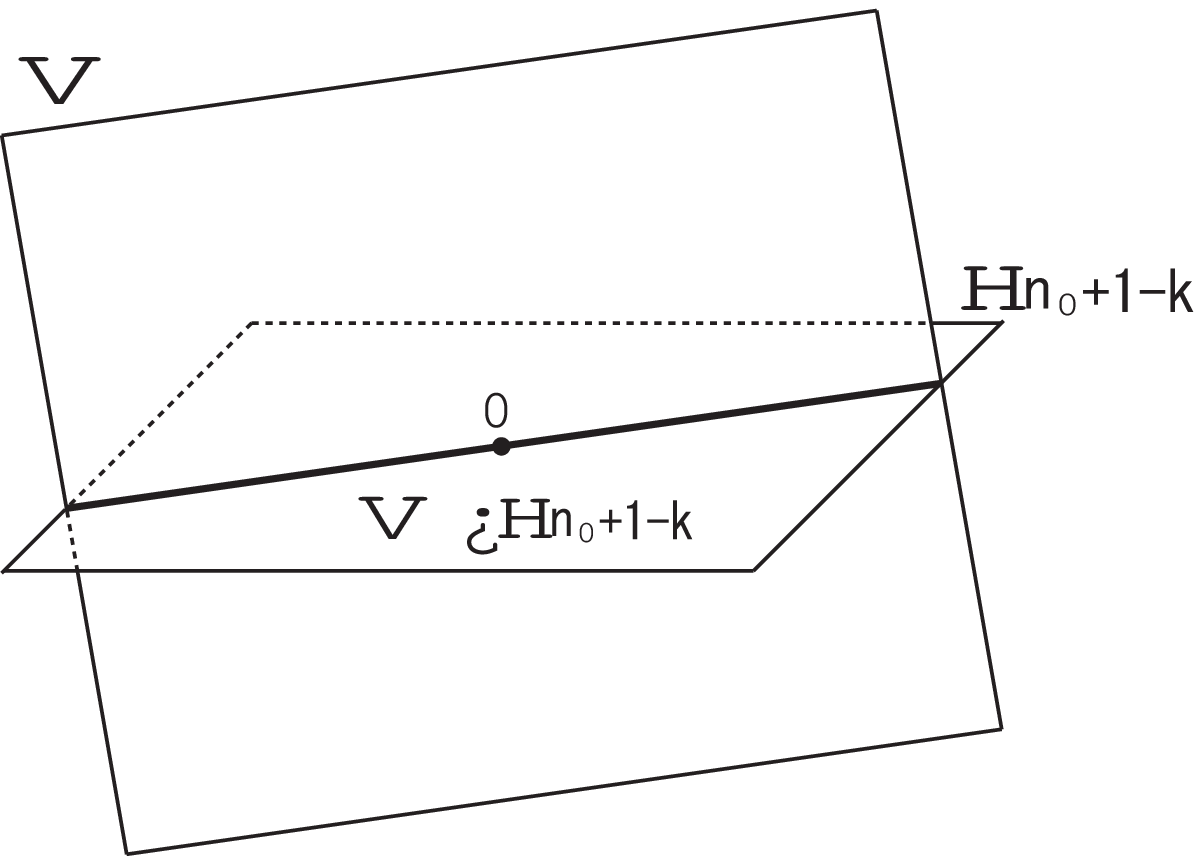}
\end{center}
\caption{}
\end{figure}

\end{lem}

Since $Gr_1 K^{n_0+1-k}=KP^{n_0-k}$ is a projective space, hence $G(n_0+1-k, K)$-paradoxical by 
Theorem \ref{thm: proj}, 
$(Gr_1 K^{n_0+1-k})^*$ is $G(n_0+1-k, K)^*$-paradoxical.
Therefore Proposition \ref{prop: key}.(5) implies that
$Gr_k K^{n_0} - (Gr_k \\K^{n_0+1-k})^*$ is $G(n_0+1-k, K)^*$-paradoxical.

On the other hand
$(Gr_k K^{n_0+1-k})^* \cong Gr_k K^{n_0+1-k}$ and 
the induction hypothesis impies $Gr_k K^{n_0+1-k}$ is $G(n_0+1-k, K)$-paradoxical,
hence $(Gr_k K^{n_0+1-k})^*$ is $G(n_0+1-k, K)^*$-paradoxical.
Therefore Proposition \ref{prop: key}.(2) implies that
$Gr_k K^{n_0}$ is $G(n_0+1-k, K)^*$-paradoxical, in particular $G(n_0, K)$-paradoxical
by Proposition \ref{prop: key}.(1).

\item
Finally we will consider the case when k satisfies $n_0/2 < k \leq n_0-1$:
Proposition \ref{prop: dual}
implies that there is a $G(n_0,K)$-equivariant homeomorphism between
$Gr_k K^{n_0}$ and $Gr_{n_0-k} K^{n_0}$,
hence Proposition \ref{prop: key}.(3) implies that
$Gr_k K^{n_0}$ is $G(n_0, K)$-paradoxical.
\end{enumerate}
\end{enumerate}

It might be an interesting question to extend the main result to
generalized flag manifolds $G/C(T)$
where $G$ is a compact and semisimple Lie group,
and $C(T)$ is a centralizer of a torus of $G$
\cite{aa}.

{\footnotesize
\renewcommand{\refname}{{\small References}}

}


\begin{thebibliography}{9}

   \bibitem[1]{aa} Andreas Arvanitoyeorgos,
   \newblock {\em An Introduction to Lie Groups and the Geometry of Homogeneous Spaces},
   \newblock Translated from the 1999 Greek original and revised by the author. Student Mathematical Library, 22. American Mathematical Society, Providence, RI, 2003.
   \bibitem[2]{bc} Richard L. Bishop and Richard J. Crittenden, 
   \newblock {\em Geometry of Manifolds},
   \newblock Pure and Applied Mathematics, Vol. XV Academic Press, New York-London 1964.
   \bibitem[3]{pg} P. G. Gormley,
   \newblock {\em Stereographic projection and the linear fractional group of transformations of quaternions},
   \newblock Proc. Roy. Irish Acad. Sect. A. 51, (1947). 67?85. (Reviewer: H. S. M. Coxeter).
   \bibitem[4]{hh} Howard Hiller, 
   \newblock {\em Geometry of Coxeter Groups},
   \newblock Research Notes in Mathematics, 54. Pitman (Advanced Publishing Program), Boston, Mass.-London, 1982.
    \bibitem[5]{sw} Stan Wagon,
   \newblock {\em The Banach-Tarski Paradox},
   \newblock With a foreword by Jan Mycielski. Corrected reprint of the 1985 original. Cambridge University Press, Cambridge, 1993.
   
\end{thebibliography}
\end{document}